\documentclass[12pt]{amsart}
\usepackage{amssymb}
\usepackage{pstcol}
\usepackage{mathptm}

\let\oldlabel=\label
\def\prellabel{\marginparsep=1em\marginparwidth=44pt
    \def\label##1{\oldlabel{##1}\ifmmode\else\ifinner\else
         \marginpar{{\footnotesize\ \\ \tt
                    ##1}}\fi\fi}}
\let\Bbb=\mathbb
\let\frak=\mathfrak

\def\a{\operatorname{a}}
\def\v{\operatorname{v}}
\def\M{{\Bbb M}}
\def\cc{{\mathfrak c}}

\def\rank{\operatorname{rank}}

\def\gp{\operatorname{gp}}
\def\Ker{\operatorname{Ker}}

\def\ini{\operatorname{in}}
\def\supp{\operatorname{supp}}

\def\Spec{\operatorname{Spec}}

\def\Aut{\operatorname{Aut}}

\def\bara{\operatorname{bar}}
\def\reg{\operatorname{reg}}

\def\Hilb{\operatorname{Hilb}}
\def\CR{\operatorname{CR}}
\def\Gr{\operatorname{Gr}}

\def\hht{\operatorname{ht}}
\def\QF{\operatorname{QF}}

\def\int{\operatorname{int}}

\let\:\colon

\def\conv{\operatorname{conv}}

\def\Q{{\Box\kern1pt}}

\let\:\colon

\let\dirsum\oplus \let\tensor\otimes

\def\RR{{\Bbb R}}

\def\QQ{{\Bbb Q}}
\def\ZZ{{\Bbb Z}}
\def\NN{{\Bbb N}}
\def\mm{{\frak m}}

\let\epsilon=\varepsilon
\let\phi=\varphi
\let\theta=\vartheta

\psset{unit=0.7cm,linewidth=0.4pt,arrowsize=1.5pt 4}
\def\vertex{\pscircle[fillstyle=solid,fillcolor=black]{0.0714}}
\newpsstyle{fatline}{linewidth=1pt}
\newpsstyle{fyp}{fillstyle=solid,fillcolor=verylight}
\definecolor{verylight}{gray}{0.95}
\definecolor{light}{gray}{0.9}
\definecolor{medium}{gray}{0.85}

\def\>#1>#2>{\xrightarrow[#2]{#1}}

\newtheorem{lemma}{Lemma}[subsection]
\newtheorem{corollary}[lemma]{Corollary}
\newtheorem{theorem}[lemma]{Theorem}
\newtheorem{proposition}[lemma]{Proposition}

\newtheorem{problem}{Problem}

\theoremstyle{definition}

\newtheorem{remark}[lemma]{Remark}

\newtheorem{example}[lemma]{Example}
\newtheorem{examples}[lemma]{Examples}

\def\subsubsection#1{\par\bigskip\noindent
\emph{#1}.\enspace}
\begin{document}
\title[Problems and algorithms for affine semigroups]
{Problems and algorithms for\\ affine semigroups}
\author{Winfried Bruns \and Joseph Gubeladze \and Ng\^o Vi\^et Trung}
\address{Universit\"at Osnabr\"uck,
FB Mathematik/Informatik, 49069 Osnabr\"uck, Germany}
\email{Winfried.Bruns@mathematik.uni-osnabrueck.de}
\address{A. Razmadze Mathematical Institute, Alexidze St. 1, 380093
Tbilisi, Georgia} \email{gubel@rmi.acnet.ge}
\address{Institute of Mathematics, Vien Toan Hoc, P. O. Box 631, Bo Ho,
Hanoi, Vietnam}
\email{nvtrung@thevinh.ncst.ac.vn\newline\rm\rule{0pt}{1.5cm}
Received December 21, 2000 and in final form July 25, 2001}

\thanks{The second author has been supported by the Deutsche
Forschungsgemeinschaft. The third author has been partially supported by the
National Basic Research Program of Vietnam.}

\dedicatory{Dedicated to the memory of Gy\"orgy
Poll\'ak\\[0.8\baselineskip]\textup{(Communicated by L\'aszl\'o M\'arki)}}

\maketitle

\section{Introduction}

Affine semigroups -- discrete analogues of convex polyhedral cones
-- mark the cross-roads of algebraic geometry, commutative algebra
and integer programming. They constitute the combinatorial
background for the theory of toric  varieties, which is their main
link to algebraic geometry. Initiated by the work of Demazure
\cite{De} and Kempf, Knudsen, Mumford and Saint-Donat \cite{KKMS}
in the early 70s, toric geometry is still a very active area of
research.

However, the last decade has clearly witnessed the extensive study
of affine semigroups from the other two perspectives. No doubt,
this is due to the tremendously increased computational power in
algebraic geometry, implemented through the theory of Gr\"obner
bases, and, of course, to modern computers.

In this article we overview those aspects of this development that
have been relevant for our own research, and pose several open
problems. Answers to these problems would contribute substantially
to the theory.

The paper treats two main topics: (1) affine semigroups and
several covering properties for them and (2) algebraic properties
for the corresponding rings (Koszul, Cohen-Macaulay, different
``sizes'' of the defining binomial ideals). We emphasize the
special case when the initial data are encoded into lattice
polytopes. The related objects -- polytopal semigroups and
algebras -- provide a link with the classical theme of
triangulations into unimodular simplices.

We have also included an algorithm for checking the semigroup
covering property in the most general setting (Section
\ref{ALGOR}). Our counterexample to certain covering conjectures
(Section \ref{CAN}) was found by the application of a small part
of this algorithm. The general algorithm could be used for a
deeper study of affine semigroups.

This paper is an expanded version of the talks given by the first
and the third author in the Problem session of the Colloquium on
Semigroups held in Szeged in July 2000.

\section{Affine and polytopal semigroups and their algebras} \label{APSA}

We use the following notation: $\ZZ$, $\QQ$, $\RR$ are the
additive groups of integral, rational, and real numbers,
respectively; $\ZZ_+$, $\QQ_+$ and $\RR_+$ denote the
corresponding additive subsemigroups of non-negative numbers, and
$\NN = \{1, 2, \dots\}$.

\subsection{Affine semigroups} An \emph{affine semigroup} is a
semigroup (always containing a neutral element) which is finitely
generated and can be embedded in $\ZZ^n$ for some $n \in \NN$.
Groups isomorphic to $\ZZ^n$ are called \emph{lattices} in the
following.

We write $\gp(S)$ for the group of differences of $S$, i.~e.\
$\gp(S)$ is the smallest group (up to isomorphism) which contains
$S$.

If $S$ is contained in the lattice $L$ as a subsemigroup, then
$x\in L$ is \emph{integral} over $S$ if $cx \in S$ for some $c \in
\NN$, and the set of all such $x$ is the \emph{integral closure}
$\bar S_L$ of $S$ in $L$. Obviously $\bar S_L$ is again a
semigroup. As we shall see in Proposition \ref{IntCl}, it is even
an affine semigroup, and can be described in geometric terms.

By a \emph{cone} in a real vector space $V=\RR^n$ we mean a subset
$C$ such that $C$ is closed under linear combinations with
non-negative real coefficients.  A cone is finitely generated if
and only if it is the intersection of finitely many vector
halfspaces. (Sometimes a set of the form $z+C$ will also be called
a cone.)  If $C$ is generated by vectors with rational or,
equivalently, integral components, then $C$ is called
\emph{rational}. This is the case if and only if the halfspaces
can be described by homogeneous linear inequalities with rational
(or integral) coefficients.

This applies especially to the cone $C(S)$ generated by $S$ in the
real vector space $L\tensor \RR$:
\begin{equation}
C(S)=\{x\in L\tensor \RR : \sigma_i(x)\ge 0,\ i=1,\dots,
s\}\tag{$*$}
\end{equation}
where the $\sigma_i$ are linear forms on $L\tensor \RR$ with
integral coefficients.

\begin{proposition}\label{IntCl}
\begin{itemize}
\item[(a)] \emph{(Gordan's lemma)}  Let $C\subset L\tensor \RR$ be a
finitely generated rational cone (i.~e.\ generated by finitely
many vectors from $L\tensor \QQ$). Then $L\cap C$ is an affine
semigroup and integrally closed in $L$.
\item[(b)]Let $S$ be an affine subsemigroup of the lattice $L$.
Then
\begin{itemize}
\item[(i)] $\bar S_L=L\cap C(S)$;
\item[(ii)] there exist $z_1,\dots,z_u\in \bar S_L$ such that $\bar
S_L=\bigcup_{i=1}^u z_i+S$;
\item[(iii)] $\bar S_L$ is an affine semigroup.
\end{itemize}
\end{itemize}
\end{proposition}

\begin{proof}
\begin{itemize}
\item[(a)] Note that $C$ is generated by finitely many elements
$x_1,\dots,x_m\in L$. Let $x\in L\cap C$. Then
$x=a_1x_1+\dots+a_mx_m$ with non-negative rational $a_i$. Set
$b_i=\lfloor a_i\rfloor$. Then
\begin{equation*}
x=(b_1x_1+\dots+b_mx_m)+(r_1x_1+\dots+r_mx_m),\qquad 0\le r_i
<1.\tag{$*$}
\end{equation*}
The second summand lies in the intersection of $L$ with a bounded
subset of $C$. Thus there are only finitely many choices for it.
These elements together with $x_1,\dots,x_m$ generate $L\cap C$.
That $L\cap C$ is integrally closed in $L$ is evident.

\item[(b)] Set $C=C(S)$, and choose a system $x_1,\dots,x_m$ of
generators of $S$. Then every $x\in L\cap C$ has a representation
$(*)$. Multiplication by a common denominator of $r_1,\dots,r_m$
shows that $x\in \bar S_L$. On the other hand, $L\cap C$ is
integrally closed by (a) so that $\bar S_L=L\cap C$.

The elements $y_1,\dots, y_u$ can now be chosen as the vectors
$r_1x_1+\dots+r_mx_m$ appearing in $(*)$. Their number is finite
since they are all integral and contained in a bounded subset of
$L\tensor\RR$. Together with $x_1,\dots,x_m$ they certainly
generate $\bar S_L$ as a semigroup.
\end{itemize}
\end{proof}

Proposition \ref{IntCl} shows that normal affine semigroups can
also be defined by finitely generated rational cones $C$: the
semigroup $S(C)=L\cap C$ is affine and integrally closed in~$L$.

We introduce special terminology in the case in which $L=\gp(S)$.
Then the integral closure $\bar S=\bar S_{\gp(S)}$ is called the
\emph{normalization}, and $S$ is \emph{normal} if $S=\bar S$.
Clearly the semigroups $S(C)$ are normal, and conversely, every
normal affine semigroup $S$ has such a representation, since
$S=S(C(S))$ (in $\gp(S)$).

Suppose that $L=\gp(S)$ and that representation $(*)$ of $C(S)$ is
irredundant. Then the linear forms $\sigma_i$ describe exactly the
support hyperplanes of $C(S)$, and are therefore uniquely
determined up to a multiple by a non-negative factor. We can
choose them to have coprime integral coefficients, and then the
$\sigma_i$ are uniquely determined. We call them the \emph{support forms}
of $S$, and write
$$
\supp(S)=\{\sigma_1,\dots,\sigma_s\}.
$$

We call a semigroup $S$ \emph{positive} if $0$ is the only
invertible element in $S$. It is easily seen that $\bar S$ is
positive as well and that positivity is equivalent to the fact
that $C(S)$ is a pointed cone with apex $0$. It is easily seen
that the map $\sigma:S\to\ZZ_+^s$,
$\sigma(x)=(\sigma_1(x),\dots,\sigma_s(x))$, is an embedding if $S$
positive. It follows that every element of $S$ can be written as
the sum of uniquely determined irreducible elements. Since $S$ is
finitely generated, the set of irreducible elements is also
finite. It constitutes the \emph{Hilbert basis} $\Hilb(S)$ of $S$;
clearly $\Hilb(S)$ is the uniquely determined minimal system of
generators of $S$. For a finitely generated positive rational cone
$C$ we set $\Hilb(C)=\Hilb(S(C))$.

Especially for normal $S$ the assumption that $S$ is positive is
not a severe restriction. It is easily seen that one has a
splitting
$$
S=S_0\dirsum S'
$$
into the maximal subgroup $S_0$ of $S$ and a positive normal
affine semigroup $S'$, namely the image of $S$ in $\gp(S)/S_0$.

\subsection{Semigroup algebras}
Now let $K$ be a field. Then we can form the \emph{semigroup
algebra} $K[S]$. Since $S$ is finitely generated as a semigroup,
$K[S]$ is finitely generated as a $K$-algebra. When an embedding
$S\to \ZZ^n$ is given, it induces an embedding $K[S]\to K[\ZZ^n]$,
and upon the choice of a basis in $\ZZ^n$, the algebra $K[\ZZ^n]$
can be identified with the Laurent polynomial ring
$K[T_1^{\pm1},\dots,T_n^{\pm1}]$. Under this identification,
$K[S]$ has the monomial basis $T^a$, $a\in S\subset \ZZ^n$ (where
we use the notation $T^a=T_1^{a_1}\cdots T_n^{a_n}$).

If we identify $S$ with the semigroup $K$-basis of $K[S]$, then
there is a conflict of notation: addition in the semigroup turns
into multiplication in the ring. The only way out would be to
avoid this identification and always use the exponential notation
as in the previous paragraph. However, this is often cumbersome.
We can only ask the reader to always pay attention to the context.

It is now clear that affine semigroup algebras are nothing but
subalgebras of $K[T_1^{\pm1},\allowbreak \dots,T_n^{\pm1}]$
generated by finitely many monomials. Nevertheless the abstract
point of view has many advantages. When we consider the elements
of $S$ as members of $K[S]$, we will usually call them
\emph{monomials}. Products $as$ with $a\in K$ and $s\in S$ are
called \emph{terms}.

The Krull dimension of $K[S]$ is given by $\rank S=\rank \gp(S)$,
since $\rank S$ is obviously the transcendence degree of
$\QF(K[S])=\QF\bigl(K[\gp(S)]\bigr)$ over $K$.

If $S$ is positive, then $\Hilb(S)$ is a minimal set of generators
for $K[S]$.

It is not difficult to check, and the reader should note that the
usage of the terms ``integral over'', ``integral closure'',
``normal'' and ``normalization'' is consistent with its use in
commutative algebra. So $K[\bar S_L]$ is the integral closure of
$K[S]$ in the quotient field $\QF(K[L])$ of $K[L]$ etc.

\subsection{Polytopal semigroup algebras}

Let $M$ be a subset of $\RR^n$. We set
\begin{align*}
L_M&=M\cap\ZZ^n,\\
E_M&=\{(x,1)\: x \in L_M\} \subset \ZZ^{n+1};
\end{align*}
so $L_M$ is the set of lattice points in $M$, and $E_M$ is the
image of $L_M$ under the embedding $\RR^n\to \RR^{n+1}$,
$x\mapsto (x,1)$. Very frequently we will consider $\RR^n$ as a
hyperplane of $\RR^{n+1}$ under this embedding; then we may
identify $L_M$ and $E_M$. By $S_M$ we denote the subsemigroup of
$\ZZ^{n+1}$ generated by $E_M$.

Now suppose that $P$ is a (finite convex) lattice polytope in
$\RR^n$, where `lattice' means that all the vertices of $P$ belong
to the integral lattice $\ZZ^n$. The affine semigroups of the type
$S_P$ will be called \emph{polytopal semigroups}. A lattice
polytope $P$ is \emph{normal} if $S_P$ is a normal semigroup.
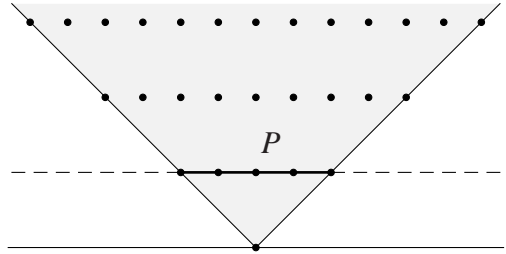
\begin{figure}[hbt]
\begin{center}
\psset{unit=1cm}
\def\vertex{\pscircle[fillstyle=solid,fillcolor=black]{0.05}}
\begin{pspicture}(-3.5,0.25)(3.5,3.25)
 \pspolygon[style=fyp,linewidth=0pt,linecolor=white]%
 (-3.25,3.25)(0,0)(3.25,3.25)
 \psline(-3.3,0)(3.3,0)
 \psline(-3.25,3.25)(0,0)(3.25,3.25)
 \psline[style=fatline](-1,1)(1,1)
 \psline[linestyle=dashed](-3.25,1)(-1,1)
 \psline[linestyle=dashed](1,1)(3.25,1)
 \rput(0,0){\vertex}
 \multirput(-1,1)(0.5,0){5}{\vertex}
 \multirput(-2,2)(0.5,0){9}{\vertex}
 \multirput(-3,3)(0.5,0){13}{\vertex}
 \rput(0.2,1.4){$P$}
\end{pspicture}
\end{center}
\caption{Vertical cross-section of a polytopal
semigroup}\label{FigSP}
\end{figure}

Let $K$ be a field. Then
$$
K[P]=K[S_P]
$$
is called a \emph{polytopal semigroup algebra} or simply a
\emph{polytopal algebra}. Since $\rank S_P=\dim P+1$ and $\dim
K[P]=\rank S_P$ as remarked above, we have
$$
\dim K[P] = \dim P + 1.
$$
Note that $S_P$ (or, more generally, $S_M$) is a graded semigroup,
i.~e.\ $S_P=\bigcup_{i=0}^\infty (S_P)_i$ such that
$(S_P)_i+(S_P)_j\subset (S_P)_{i+j}$; its $i$-th graded component
$(S_P)_i$ consists of all the elements $(x,i)\in S_P$. Moreover,
$S_P$ is even \emph{homogeneous}, namely generated by its elements
of degree $1$.

Therefore $R=K[P]$ is a graded $K$-algebra in a natural way and
generated by its degree $1$ elements. Its $i$-th graded component
$R_i$ is the $K$-vector space generated by $(S_P)_i$. The elements
of $E_P=(S_P)_1$ have degree $1$, and therefore $R$ is a
homogeneous $K$-algebra in the terminology of Bruns and Herzog
\cite{BH}. The defining relations of $K[P]$ are the binomials
representing the affine dependencies of the lattice points of $P$.
(In Section \ref{Trung} we will discuss the properties of the
ideal generated by the defining binomials.) Some easy examples:

\begin{examples}
\begin{itemize}
\item[(a)] $P=\conv(1,4)\in\RR^1$. Then $P$ contains the lattice
points $1, 2, 3,\allowbreak 4$, and the relations of the corresponding
generators of $K[P]$ are given by
$$
X_1X_3=X_2^2,\ X_1X_4=X_2X_3,\ X_2X_4=X_3^2.
$$
\item[(b)] $P=\conv\big((0,0), (0,1), (1,0), (1,1)\big)$. The lattice
points of $P$ are exactly the 4 vertices, and the defining
relation of $K[P]$ is $X_1X_4=X_2X_3$.
\item[(c)] $P=\conv\big((1,0), (0,1), (-1,-1)\big)$. There is a fourth
lattice point in $P$, namely $(0,0)$, and the defining relation is
$X_1X_2X_3=Y^3$ (in suitable notation).
\end{itemize}
\end{examples}

\begin{figure}[htb]
\begin{center}
\begin{pspicture}(0,-0.5)(3,0.5)
\multirput(0,0)(1,0){4}{\vertex} \psline(0,0)(3,0)
\end{pspicture}\qquad\qquad
\begin{pspicture}(0,0)(1,1)
\pspolygon[style=fyp](0,0)(1,0)(1,1)(0,1)
\multirput(0,0)(1,0){2}{\vertex} \multirput(0,1)(1,0){2}{\vertex}
\end{pspicture}\qquad\qquad
\begin{pspicture}(-1,-0.5)(1,1.5)
\pspolygon[style=fyp](-1,-1)(1,0)(0,1) \rput(-1,-1){\vertex}
\rput(0,1){\vertex} \rput(1,0){\vertex} \rput(0,0){\vertex}
\end{pspicture}
\end{center}
\caption{}\label{ThreePol}
\end{figure}

Note that the polynomial ring $K[X_1,\dots,X_n]$ is a polytopal
algebra, namely $K[\Delta_{n-1}]$ where $\Delta_{n-1}$ denotes the
$(n-1)$-dimensional unit simplex.

It is often useful to replace a polytope $P$ by a multiple $cP$
with $c\in \NN$. The lattice points in $cP$ can be identified with
the lattice points of degree $c$ in the cone $C(S_P)$; in fact,
the latter are exactly of the form $(x,c)$ where $x\in L_{cP}$.

Polytopal semigroup algebras appear as the coordinate rings of
projective toric varieties; see Oda \cite{Oda}

\section{Hilbert bases of affine normal semigroups}\label{CAN}

\subsection{Normality and covering}\label{CANInt}

In this section we will investigate the question whether the
normality of a positive affine semigroup can be characterized in
terms of combinatorial conditions on its Hilbert basis.

Let $C$ be a cone in $\RR^n$ generated by finitely many rational
(or integral) vectors. We say that a collection of rational
subcones $C_1,\dots,C_m$ is a \emph{triangulation} of $C$ if $C_i$
is simplicial for all $i$ (i.e.\ generated by a linearly
independent set of vectors), $C=C_1\cup\dots\cup C_m$ and
$C_{i_1}\cap\dots\cap C_{i_k}$ is a face of $C_{i_1},\dots,C_{i_k}$
for every subset $\{i_1,\dots,i_k\}\subset\{1,\dots,m\}$.

Let $M$ be a subset of a cone $C$ as above. An
\emph{$M$-triangulation of $C$} is a triangulation into simplicial
cones spanned by subsets of $M$, and a \emph{Hilbert
triangulation} is a $\Hilb(S(C))$-triangulation of $C$.

Correspondingly, a \emph{Hilbert subsemigroup} $S'$ of $S$ is a
subsemigroup generated by a subset of $\Hilb(S)$. We say that $S$
is \emph{covered} by subsemigroups $S_1,\dots, S_m$ if
$S=S_1\cup\dots\cup S_m$.

A subset $X$ of $\ZZ^n$ is called \emph{unimodular} if it is part
of a basis of $\ZZ^n$; in other words, if it is linearly
independent and generates a direct summand of $\ZZ^n$. Cones and
semigroups are \emph{unimodular} if they are generated by
unimodular sets, and a collection of unimodular objects is
likewise called unimodular.

\begin{proposition}\label{UHCsat}
If $S$ is covered by unimodular subsemigroups, then it is normal.
More generally, if $S$ is the union of normal subsemigroups $S_i$
such that $\gp(S_i)=\gp(S)$, then $S$ is also normal.
\end{proposition}

This follows immediately from the definition of normality.

We will see in Corollary \ref{uh=ic} that the hypothesis
$\gp(S_i)=\gp(S)$ is superfluous, and that it is even enough that
the $S_i$ cover $S$ ``asymptotically'' .

The following converse is important for the geometry of toric
varieties; it provides the combinatorial basis for the equivariant
resolution of their singularities.

\begin{theorem}\label{triang}
Every finitely generated rational cone $C\subset \RR^n$ has a
unimodular triangulation.
\end{theorem}

It is not difficult to prove the theorem for which we may assume
that $\dim C=n$. One starts with an arbitrary triangulation of
$C$, and considers each of the involved simplicial subcones $C'$.
The shortest nonzero integer vectors on each of the rays of $C'$ form a
linearly independent set $X$. If $X$ is not unimodular, then $X$
is not the Hilbert basis of $S(C')$, and one subdivides $C'$ by
one of the vectors $r_1x_1+\dots+r_mx_m$ appearing in the proof of
Gordan's lemma. For each of the simplicial subcones $C''$
generated by subdivision the group $\gp(S(C''))$ has smaller index
than $\gp(S(C'))$ in $\ZZ^n$.  After finitely many steps one thus
arrives at a unimodular triangulation.

Especially for polytopal semigroups, Theorem \ref{triang} is not
really satisfactory, since it is not possible to interpret it in
the lattice structure of a polytope $P\subset \ZZ^n$. In fact,
only the simplicial Hilbert subcones of $C(S_P)$ correspond to the
lattice simplices contained in $P$. It is not hard to see that the
cone spanned by a lattice simplex $\delta\subset P$ is unimodular
if and only if $\delta$ has the smallest possible volume $1/n!$.
Such simplices are also called \emph{unimodular}. Furthermore,
$P$ (regardless of its dimension) can be triangulated into \emph{empty}
lattice simplices, i.~e.\ simplices $\delta$ such that $\delta\cap\ZZ^n$
is exactly the set of vertices of $\delta$.

Suppose now that $P$ is a lattice polytope of dimension $2$ and
triangulate it into empty lattice simplices. Since, by Pick's
theorem, an empty simplex of dimension $2$ has area $1/2$, one
automatically has a unimodular triangulation. It follows
immediately that $S_P$ is the union of unimodular Hilbert
subsemigroups and thus normal. Moreover, $C(S_P)$ has a unimodular
Hilbert triangulation.
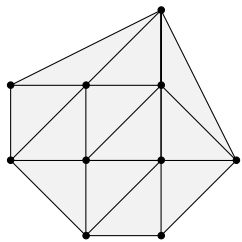
\begin{figure}[hbt]
$$
\psset{unit=1cm}
\def\vertex{\pscircle[fillstyle=solid,fillcolor=black]{0.05}}
\begin{pspicture}(-0.5,.25)(3.0,2,75)
 \pspolygon[style=fyp](1,0)(2,0)(3,1)(2,3)(0,2)(0,1)
 \multirput(1,0)(1,0){2}{\vertex}
 \multirput(0,1)(1,0){4}{\vertex}
 \multirput(0,2)(1,0){3}{\vertex}
 \multirput(2,3)(1,0){1}{\vertex}
 \psset{fillstyle=none}
 \psline(1,1)(2,1)(2,2)(2,3)(2,2)(1,2)(1,1)
 \psline(2,3)(0,1)(1,1)(1,0)(2,1)(3,1)(2,2)
 \psline(2,0)(2,1)
 \psline(0,2)(1,2)
 \psline(1,1)(2,1)(2,2)(1,1)
\end{pspicture}
$$
\caption{Triangulation of a lattice polygon}
\end{figure}

More generally, Seb\H o has shown the following

\begin{theorem}\label{UHT3}
Every positive finitely generated cone of dimension $3$
has a unimodular Hilbert triangulation.
\end{theorem}

We refer the reader to Seb\H os paper \cite{Se} or to \cite{BGgren}
for the proof, which is by no means straightforward. The much
simpler polytopal case discussed above is characterized by the
fact that the elements of the Hilbert basis of $C(S)$ lie in a
hyperplane.

Theorem \ref{UHT3} also holds in dimension $1$ and $2$ where it is
easily proved, but it cannot be extended to dimension $\ge 4$, as
shown by a counterexample due to Bouvier and Gonzalez-Sprinberg
\cite{BoGo}.

As has been mentioned already, triangulations are very interesting
objects for the geometry of toric varieties. Triangulations also
provide the connection between discrete geometry and Gr\"obner
bases of the binomial ideal defining a semigroup algebra. See
Sturmfels \cite{Stu} for this important and interesting theme; we
will briefly discuss it in Section \ref{Trung}.

Despite of counterexamples to the existence of unimodular Hilbert
triangulations in dimension $\ge 4$, it is still reasonable to
consider the following, very natural sufficient condition of \emph{unimodular
Hilbert covering} for positive normal semigroups $S$:
\medskip

\noindent ({\bf UHC})\enspace \emph{$S$ is covered by its
unimodular Hilbert subsemigroups.}
\medskip

For polytopal semigroups (UHC) has a clear geometric
interpretation: it just says that $P$ is the union of its
unimodular lattice subsimplices.

Seb\H o \cite[Conjecture B]{Se} has conjectured that (UHC) is
satisfied by all normal affine semigroups. Below we present a
$6$-dimensional counterexample to Seb\H o's conjecture. However
there are also positive results on (UHC) and even on unimodular
triangulations for multiples $cP$ of polytopes; see Subsection
\ref{HighMult}.

A natural variant of (UHC), and weaker than (UHC), is the
existence of a \emph{free Hilbert cover}:
\medskip

\noindent ({\bf FHC})\enspace \emph{$S$ is the union (or covered
by) the subsemigroups generated by the linearly independent
subsets of $\Hilb(S)$.}
\medskip

For (FHC) -- in contrast to (UHC) --  it is not evident that it
implies the normality of the semigroup. Nevertheless it does so,
as we will see in Corollary \ref{uh=ic}. A formally weaker -- and
certainly the most elementary -- property is the \emph{integral
Carath\'eodory property}:
\medskip

\noindent ({\bf ICP})\enspace \emph{Every element of $S$ has a
representation $x=a_1s_1+\dots+a_ms_m$ with $a_i\in\ZZ_+$,
$s_i\in\Hilb(C)$, and $m\le \rank S$. }\medskip

Here we have borrowed the well-motivated terminology of Firla and
Ziegler \cite{FZ}: (ICP) is obviously a discrete variant of
Carath\'eodory's theorem for convex cones. It was first asked in
Cook, Fonlupt, and Schrijver \cite{CFS} whether all cones have
(ICP) and then conjectured in \cite[Conjecture A]{Se} that the
answer is `yes'.

Later on we will use the \emph{representation length}
$$
\rho(x)=\min\{m \mid x=a_1s_1+\dots+a_ms_m,\ a_i\in\ZZ_+,\ s_i\in
\Hilb(S)\}
$$
for an element $x$ of a positive affine semigroup $S$. If $\rho(x)\le m$,
we also say that $x$ is $m$-\emph{represented}. In order to measure
the deviation of $S$ from (ICP), we introduce the notion of
\emph{Carath\'eo\-dory rank\/} of an affine semigroup $S$,
$$
\CR(S)=\max\{\rho(x) \mid x\in S\}.
$$

Variants of this notion, called asymptotic and virtual
Carath\'eodory rank will be introduced in Section \ref{ALGOR}.

The following $10$ vectors constitute the Hilbert basis of a
normal positive semigroup $S_6$:
\begin{alignat*}{2}
 z_1 & = (0,\,1,\,0,\,0,\,0,\,0), & \qquad z_6 &=
(1,\,0,\,2,\,1,\,1,\,2), \\
 z_2 & = (0,\,0,\,1,\,0,\,0,\,0), & \qquad z_7 &=
(1,\,2,\,0,\,2,\,1,\,1), \\
 z_3 & = (0,\,0,\,0,\,1,\,0,\,0), & \qquad z_8 &=
(1,\,1,\,2,\,0,\,2,\,1), \\
 z_4 & = (0,\,0,\,0,\,0,\,1,\,0), & \qquad z_9 &=
(1,\,1,\,1,\,2,\,0,\,2), \\
 z_5 & = (0,\,0,\,0,\,0,\,0,\,1), & \qquad z_{10} &=
(1,\,2,\,1,\,1,\,2,\,0).
\end{alignat*}
As a counterexample to (UHC) it was found by the first two authors
\cite{BGcov}. In cooperation with Henk, Martin and Weismantel
\cite{BGHMW} it was then shown that $\CR(S_6)=7$ so that (ICP)
does not hold for all normal affine semigroups $S$. The cone $C_6$
and the semigroup $S_6=S(C_6)$ have several remarkable properties;
for example, $\Aut(S_6)$ operates transitively on the Hilbert
basis. The reader can easily check that $z_1,\dots,z_{10}$ lie on
a hyperplane. Therefore $S_6=S_P$ for a $5$-dimensional lattice
polytope $P$. Further details can be found in the papers just
quoted.

A crucial idea in finding $S_6$ was the introduction of the class
of \emph{tight} cones and semigroups; see \cite{BGcov}.

So far one does not know a semigroup $S$ satisfying (ICP), but not
(UHC). This suggests the following problem:

\begin{problem}
Does \emph{(ICP)} imply \emph{(UHC)}?
\end{problem}

Since the positive results end in dimension $3$ and the
counterexample lives in dimension $6$, the situation is completely
open in dimensions $4$ and $5$:

\begin{problem}
Prove or disprove \emph{(ICP)} and/or \emph{(UHC)} in dimension $4$
and $5$.
\end{problem}

We have seen above that every triangulation of a lattice polygon
into empty lattice simplices is unimodular. This property is truly
restricted to dimension at most $2$. In fact, Hosten, MacLagan,
and Sturmfels \cite{HMS} have given an example of a
$3$-dimensional cone that contains no finite set $M$ of lattice
points such that every triangulation of $C$ using all the points
of $M$ is unimodular.

\subsection{An upper bound for Carath\'eodory rank}\label{UpperCR}

Let $p_1,\dots,p_n$ be different prime numbers, and set
$q_j=\prod_{i\neq j} p_i$. Let $S$ be the subsemigroup of $\ZZ_+$
generated by $q_1,\dots,q_n$. Since $\gcd(q_1,\dots,q_n)=1$, there
exists an $m\in\ZZ_+$ with $u\in S$ for all $u\ge m$. Choose $u\ge
m$ such that $u$ is not divisible by $p_i$, $i=1\dots,n$. Then all
the $q_i$ must be involved in the representation of $u$ by
elements of $\Hilb(S)$. This example shows that there is no bound
of $\CR(S)$ in terms of $\rank S$ without further conditions on
$S$.

For normal $S$ there is a linear bound for $\CR(S)$ as given by
Seb\H o \cite{Se}:

\begin{theorem}\label{SebCR}
Let $S$ be a normal positive affine semigroup of rank $\ge 2$.
Then $\CR(S)\le 2(\rank(S)-1)$.
\end{theorem}

\begin{figure}[hbt]
$$
\begin{pspicture}(-0.5,0.0)(6.3,6.8)
 \pspolygon[linecolor=white,fillcolor=medium,fillstyle=solid,linewidth=0pt]%
 (2.5,6.25)(2,5)(1,2)(1,1)(2,1)(6.3,3.15)(6.3,6.35)
 \pspolygon[fillcolor=verylight,fillstyle=solid,linecolor=white,linewidth=0pt]%
 (2,5)(1,2)(1,1)(2,1)(0,0)
 \multirput(1,1)(1,0){2}{\vertex}
 \multirput(1,2)(1,0){4}{\vertex}
 \multirput(2,3)(1,0){4}{\vertex}
 \multirput(2,4)(1,0){5}{\vertex}
 \multirput(2,5)(1,0){5}{\vertex}
 \multirput(3,6)(1,0){4}{\vertex}
 \rput(0,0){\vertex}
 \psset{fillstyle=none}
 \psline[linecolor=darkgray](0,0)(6.3,3.15)
 \psline[linecolor=darkgray](0,0)(2.5,6.25)
 \psline[style=fatline](2,5)(1,2)(1,1)(2,1)
 \rput(3.5,3.5){$C'(S)$}
\end{pspicture}
$$
\caption{The bottom} \label{FigBot}
\end{figure}
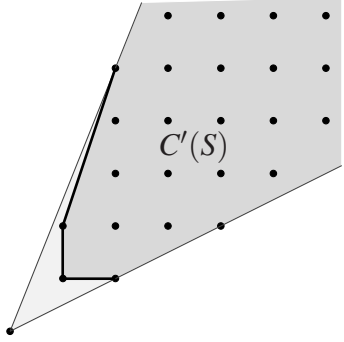

For the proof we denote by $C'(S)$ the convex hull of
$S\setminus\{0\}$ (in $\gp(S)\tensor\RR$). Then we define the
\emph{bottom} $B(S)$ of $C'(S)$ by
$$
B(S)=\bigl\{x\in C'(S): [0,x]\cap C'(S)=\{x\}\bigr\}
$$
($[0,x]=\conv(0,x)$ is the line segment joining $0$ and $x$). In
other words, the bottom is exactly the set of points of $C'(S)$
that are \emph{visible} from $0$ (see Figure \ref{FigBot}).

Let $H$ be a support hyperplane intersecting $C'(S)$ in a compact
facet. Then there exists a unique primitive $\ZZ$-linear form
$\gamma$ on $\gp(S)$ such that $\gamma(x)=a>0$ for all $x\in H$
(after the extension of $\gamma$ to $\gp(S)\tensor \RR$). Since
$\Hilb(S)\cap H\neq\emptyset$, one has $a\in \ZZ$. We call
$\gamma$ the \emph{basic grading} of $S$ associated with the facet
$H\cap C'(S)$ of $C'(S)$. It can be thought of as the graded
structure
$$
\deg_\gamma:S\to\ZZ_+,\qquad x\mapsto\gamma(x).
$$
\begin{proof}[Proof of Theorem \ref{SebCR}] It is easily seen that
the bottom of $S$ is the union of finitely many lattice polytopes
$F$, all of whose lattice points belong to $\Hilb(S)$. We now
triangulate each $F$ into empty lattice subsimplices. Choose $x\in
S$, and consider the line segment $[0,x]$. It intersects the
bottom of $S$ in a point $y$ belonging to some simplex $\sigma$
appearing in the triangulation of a compact facet $F$ of $C'(S)$.
Let $z_1,\dots,z_n\in \Hilb(S)$, $n=\rank(S)$, be the vertices of
$\sigma$. Then we have
$$
x=(a_1z_1+\dots+a_nz_n)+(q_1z_1+\dots+q_nz_n),\qquad a_i\in
\ZZ_+,\ q_i\in\QQ,\ 0\le q_i<1,
$$
as in the proof of Gordan's lemma. Set $x'=\sum_{i=1}^n q_iz_i$,
let $\gamma$ be the basic grading of $S$ associated with $F$, and
$a=\gamma(y)$ for $y\in F$. Then $\gamma(x')<na$, and at most
$n-1$ elements of $\Hilb(S)$ can appear in a representation of
$x'$. This shows that $\CR(S)\le 2n-1$.

However, this bound can be improved. Set $x''=x_1+\dots+x_n-x'$.
Then $x''\in S$, and it even belongs to the cone generated by
$x_1,\dots,x_n$. If $\gamma(x'')< a$, one has $x''=0$. If
$\gamma(x'')=a$, then $x''$ is a lattice point of $\sigma$. By the
choice of the triangulation this is only possible if $x''=x_i$ for
some $i$, a contradiction. Therefore $\gamma(x'')>a$, and so
$\gamma(x')<(n-1)a$. It follows that $\CR(S)\le 2n-2$.
\end{proof}

In view of Theorem \ref{SebCR} it makes sense to set
$$
\mathcal{CR}(n)=\max\bigl\{\CR(S): S \text{ is normal positive and
}\rank S=n\bigr\}.
$$
With this notion we can reformulate Theorem \ref{SebCR} as
$\mathcal{CR}(n)\le 2(n-1)$. On the other hand, the
counterexample $S_6$ to (ICP) presented above implies that
$$
\mathcal{CR}(n)\ge \left\lfloor \frac {7n}6 \right\rfloor.
$$
In fact, $\rank S_6=6$ and $\CR(S_6)=7$. Therefore suitable direct
sums $S_6\oplus\dots\oplus S_6\dirsum\ZZ_+^p$ attain the lower
bound just stated.

\begin{problem}
Improve one or both of the inequalities for $\mathcal{CR}(n)$.
\end{problem}

\subsection{Unimodular covering of high multiples of polytopes}
\label{HighMult}

The counterexample presented above shows that a normal lattice
polytope need not be covered by its unimodular lattice
subsimplices. However, this always holds for a sufficiently high
multiple of $P$ \cite{BGT}:

\begin{theorem}\label{A1.5}
For every lattice polytope $P$ there exists $c_0 > 0$ such that
$cP$ is covered by its unimodular lattice subsimplices (and,
hence, is normal by Proposition \ref{UHCsat}) for all $c \in \NN$,
$c > c_0$.
\end{theorem}

A proof can be found in \cite{BGT} or \cite{BGgren}. For
elementary reasons one can take $c=1$ in dimension $1$ and $2$,
and it was communicated by Ziegler that $c=2$ suffices in
dimension $3$. This is proved by Kantor and Sarkaria \cite{KS};
moreover, they show that $4P$ has a unimodular triangulation for
every lattice polytope $P$ in dimension $3$.

\begin{problem}
Is it possible to choose $c_0$ only depending on the dimension of
$P$? If the answer is positive, give an explicit estimate for
$c_0$ in terms of $\dim P$.\footnotemark\footnotetext{Problem 4 has 
meanwhile been solved positively. See W. Bruns and J. Gubeladze, 
{\em Unimodular covers of multiples of polytopes} (in preparation), where a 
subexponential bound for $c_0$ is given.}
\end{problem}

For normality (without unimodular covering) this problem has a
satisfactory answer:

\begin{theorem}\label{cPnorm}
For every lattice polytope $P$ the multiples $cP$ are normal for
$c\ge \dim P-1$.
\end{theorem}

This can be shown by essentially the same arguments as Theorem
\ref{SebCR}; see \cite{BGT} for another argument.

In fact, it is proved in \cite{KKMS} that one even has a stronger
statement on the existence of unimodular Hilbert triangulations:

\begin{theorem}\label{cPtri}
For every lattice polytope $P$ there exists $c_0 > 0$ such that
$cP$ has a unimodular triangulation for all multiples $c=kc_0$,
$k\in\NN$.
\end{theorem}

However, note that Theorem \ref{A1.5} makes an assertion on all
sufficiently large $c$, whereas Theorem \ref{cPtri} only concerns
the multiples of a single $c_0>0$:

\begin{problem}
Does $cP$ have a unimodular triangulation for all $c\gg 0$?
\end{problem}

For applications in algebraic geometry or commutative algebra one
is especially interested in so-called \emph{regular} (or
projective) triangulations. We will come back to this point in
Section \ref{Trung}.

\section{Algorithms for coverings}\label{ALGOR}

An affine semigroup $S$ is a subset of a free abelian group
equipped with a minimal amount of algebraic structure, but this
suffices to specify $S$ by finite data, namely a generating set.
Therefore, the question of deciding whether an affine semigroup is
the union of a given system of sub-semigroups, also specified in
terms of generators, seems interesting. In this section we develop
an algorithm deciding in a finite number of steps whether $S$ is
covered by a system of subsemigroups. Actually, in the process of
checking this property we have to treat the more general situation
of ``modules'' over affine semigroups. The connection with
Carath\'eodory ranks and (ICP) will also be outlined.

The algorithm contains subalgorithms for checking
\emph{asymptotic} and \emph{virtual} covering properties.

For subsets $A,B\subset\ZZ^n$ we use the following notation
$$
\pi(A|B)=\lim_{\epsilon\to\infty}\frac{\#\{a\in A\cap B:\
\|a\|<\epsilon\}} {\#\{b\in B:\ \|b\|<\epsilon\}}
$$
provided the limit exists. (Here $\|-\|$ denotes the standard
Euclidean norm in $\RR^n$.) One should interpret $\pi(A|B)$ as the
probability with which a random element of $B$ belongs to $A$.

 From the view point of geometry it is preferable to associate
objects in $\RR^n$ with polytopes and cones. However, the reader
should note that all data are specified in terms of rational
vectors, and that the algorithms below only require arithmetic
over $\QQ$ (or $\ZZ$).

\subsection{Normal affine semigroups}

For the algorithms developed below it is important that certain
basic computations for normal semigroups can be carried out:
\begin{itemize}
\item[(a)] The determination of the Hilbert basis of $S(C)$ where $C$
is the cone given by finitely many elements
$z_1,\dots,z_m\in\ZZ^n$. They generate the integral closure of the
affine semigroup $\ZZ_+z_1+\dots+\ZZ_+z_m$.

\item[(a')] The determination of a finite system of generators of
$S(C)$ as a module over the semigroup generated by
$z_1,\dots,z_m$.

\item[(b)] The description of the cone $C$ by a system of homogeneous
rational inequalities.

\item[(c)] The reverse process of determining $\Hilb(C)$ from a
description of $C$ by inequalities.

\item[(d)] The computation of a triangulation of $C$ into
simplicial subcones spanned by elements of $\{z_1,\dots,z_m\}$.
\end{itemize}

Note that the computations (b) and (c) are dual to each other
under exchanging $C$ with its dual cone
$C^*=\{\phi\in(\RR^n)^*:\phi(C)\subset\RR_+\}$. Nevertheless one
should mention (c) explicitly, since it allows one to compute
intersections of cones.

Algorithms for (a)--(d) have already been implemented in NORMALIZ
\cite{BK}, and the documentation of this program describes the
details. In the following we will refer to NORMALIZ whenever one
of these computations has to be carried out.

\subsection{Asymptotic covers}\label{ASYMP}

Let $S\subset\ZZ^n$ be an affine semigroup, neither necessarily
positive nor necessarily of full rank $n$. A subset
$M\subset\ZZ^n$ is called an \emph{$S$-module} if $S+M\subset M$.
A module $M$ is called \emph{finitely generated} if
$M=\{m_1+s,\dots,m_k+s: s\in S\}$ for some finite subset
$\{m_1,\dots,m_k\}\subset M$. For finitely generated modules we
write $M\in\M(S)$. Notice, that in the special case $S=0$ a
finitely generated $S$-module is just a finite set (maybe
$\emptyset$).

Consider an affine semigroup $S$ and a finite family of affine
semigroups
$$
S_1,\dots,S_t\subset S.
$$
We say that  $S$ is covered
\emph{asymptotically} by the $S_i$ if $\pi(S_1\cup\cdots\cup
S_t|S)=1$. One should observe that the notion of asymptotic covering
is an intrinsic property of the semigroup $S$ and the family
$\{S_1,\dots,S_t\}$. In other words, it does not depend on the
embedding $S\to\ZZ^n$. Further, $S$ is said to be \emph{virtually}
covered by the $S_i$ if $\#\bigl(S\setminus(S_1\cup\dots\cup
S_t)\bigr)<\infty$.

Now assume we are given a finitely generated $S$-module $M$ and
$S_i$-submodules $M_i\subset M$ so that $M_i\in\M(S_i)$
$i\in[1,t]$. One then introduces the notions of covering,
asymptotic covering and virtual covering of $M$ by the $M_i$ in
the obvious way.

\begin{lemma}\label{cond}
For an affine semigroup $S$ the conductor ideal $\cc_{\bar S/S}=
\{x\in S:\ x+\bar S\subset S\}$ is a nonempty set.
\end{lemma}

\begin{proof}
Let $G$ be a generating set of $S$ and $\bar G$ be a finite
generating set of $\bar S$ as a module over $S$. That $\bar S$ is
in fact a finitely generated $S$-module, has been stated in Lemma
\ref{IntCl}. Fix representations $z=x_z-y_z$, $z\in\bar G$,
$x_z,y_z\in G$. Then $\sum_{z\in G}y_z\in \cc_{\bar S/S}$.
\end{proof}

Since one can effectively compute a system of generators of the
$S$-module  $\bar S$ once a generating set of $S$ is given, the
proof of Lemma \ref{cond} provides an algorithm for computing an
element of $\cc_{\bar S/S}$ if a generating set of $S$ is given.
This algorithm is called CONDUCTOR.

Consider an affine semigroup $S\subset\ZZ^n$ and a family of
affine sub-semigroups $S_1,\dots,S_t\subset S$, $t\in\NN$. Their
cones in $\RR^n$ will be denoted correspondingly by
$C(S),C(S_1),\dots,C(S_t)$. A family of \emph{non-empty} modules
$M\in\M(S)$, $M_1\in\M(S_1),\dots,M_t\in\M(S_t)$, such that
$M_1,\dots,M_t\subset M$ ($\subset\ZZ^n$), is also assumed to be
given.

Put
$$
\Sigma=\biggl\{\sigma\subset[1,t]:\
\dim(\bigcap_{i\in\sigma}C(S_i))=\rank S\ \ \text{and}\ \
\bigcup_{i\in\sigma}\gp(S_i)=\gp(S)\biggr\}
$$
and
$$
C_\sigma=\bigcap_{i\in\sigma}C(S_i),\qquad\sigma\in\Sigma.
$$

\begin{lemma}\label{critas}
$S$ is asymptotically covered by $S_1,\dots,S_t$ if and only if
$C(S)=\bigcup_\Sigma C_\sigma$. Moreover, $M$ is asymptotically
covered by the $M_i$ if and only if the following implication holds
for every $z\in\ZZ^n$:
\begin{multline*}
(z+\gp(S))\cap M\neq\emptyset)\quad\implies\\
 S\ \text{is asymptotically covered by}\ \{S_j:\ j\in[1,t],\
    (z+\gp(S))\cap M_j\neq\emptyset\}.
\end{multline*}
\end{lemma}

\begin{proof}
Consider finite generating sets $G_i\subset S_i$, $i\in[1,n]$. The
affine hyperplanes in $\RR\tensor\gp(S)$, spanned by the elements
of $\bigcup_1^nG_i$, cut the cone  $C(S)$ into subcones which we
call \emph{elementary cells}, i.~e. the elementary cells are the
maximal dimensional cones in the obtained polyhedral subdivision
of $C(S)$. Clearly, the elementary cells are again finite rational
cones. So by Gordan's lemma the semigroups $S\cap C$ are all
affine. (The general form of Gordan's lemma used here and below follows from
\ref{IntCl} and \cite[7.2]{BGdiv}.)

$S$ is asymptotically covered by the $S_i$ if and only if
$\pi\bigl(S_1\cup\dots\cup S_n\big|S\cap E \bigr)=1$ for every
elementary cell $E$, or equivalently
$$
\pi\biggl(\bigcup_{i\in\sigma_E}S_i\cap E\bigg|S\cap E\biggr)=1
$$
where $\sigma_E=\{i\in[1,n]:\ E\subset C(S_i)\}$, $E$ running
through the set elementary cells.

We claim that $S$ is asymptotically covered if and only if
$\sigma_E\in\Sigma$. This clearly proves the first part of the
lemma.

The ``only if'' part of the claim follows easily from the fact that
$\gp(S\cap E)=\gp(S)$. For the ``if'' part we pick elements
$z_i\in\cc_{\bar S_i/S_i}$, $i\in\sigma_E$ (Lemma \ref{cond}).
Then the assumption $\sigma_E\in\Sigma$ implies
$$
S_0:=\gp(S)\cap
E\cap\biggl(\bigcap_{i\in\sigma_E}\bigl(z_i+C(S_i)\bigr)\biggr)\subset
S\cap E
$$
and we are done because by elementary geometric consideration one
has
$$
\pi\biggl(\bigcup_{i\in\sigma_E}S_i\cap E \bigg|S_0\biggr)=1.
$$
Now assume the implication $\implies$ of the lemma holds. $M$ is
contained in finitely many residue classes in $\ZZ^n$ modulo
$\gp(S)$. By fixing origins in these classes and taking
intersections with the modules $M,M_1,\dots,M_t$, the general case
reduces to the situation when $M,M_1,\dots,M_t\subset\gp(S)$. Pick
elements $y_i\in M_i$. Then we have
$$
M_\sigma:=\gp(S)\cap\bigcap_{i\in\sigma}\bigl(y_i+z_i+C(S_i)\bigr)\subset
M,\qquad \sigma\in\Sigma,
$$
with the $z_i$ chosen as above. We are done by the following
observations:
$$
M_\sigma\subset\bigcup_{i\in\sigma}M_i
$$
and
$$
\pi\biggl(\bigcup_{\sigma\in\Sigma}M_\sigma\bigg| M\biggr)=1,
$$
the latter equality being easily deduced from the condition
$C(S)=\bigcup_\Sigma C_\sigma$.

Now assume $M$ is asymptotically covered by the $M_i$. Then
we have the implication
\begin{multline*}
(z+\gp(S))\cap M\neq\emptyset)\quad\implies\quad
(z+\gp(S))\cap M\ \text{is asymptotically covered by}\\
\{(z+\gp(S))\cap M_j:\ j\in[1,t],\
(z+\gp(S))\cap M_j\neq\emptyset\}.
\end{multline*}
It only remains to notice that each of these $(s+\gp(S))\cap M_j$
is asymptotically covered by $m_j+S$ for an arbitrary element
$m_j\in(z+\gp(S))\cap M_j$, and, similarly, $(z+\gp(S))\cap M$ is
asymptotically covered by $m+S$, $m\in(z+\gp(S))\cap M$.
\end{proof}

The proof of Lemma \ref{critas} gives an algorithm deciding
whether $S$ is asymptotically covered by $S_1,\dots S_t$, using
explicit generating sets as input. In fact, the conditions (i)
that a finite rational cone is covered by a system of finite
rational subcones and (ii) that a finitely generated free abelian
group is covered by a system of subgroups, can both be checked
effectively. It is of course necessary that we are able to compute
the cone of an affine semigroup once a generating set of the
semigroup is given (NORMALIZ), to form the intersection of a
system of finite rational cones (given in terms of the support
inequalities) and, furthermore, to compute the group of differences
of an affine semigroups.

In fact, for the cone covering property we first triangulate the
given cone $C$ (using only extreme generators) and then inspect
successively the resulting simplicial subcones as follows. If such
a simplicial cone $T$ is contained in one of the given cones, say
$C_1,\dots,C_t$, it is neglected and we pass to another simplicial
cone. If it is not contained in any of the cones $C_1,\dots,C_t$,
then we split $T$ into two cones (of the same dimension) by the
affine hull of a facet $F\subset C_i$ for some $i\in[1,t]$.
Thereafter the two pieces of $T$ are tested for the containment
property in one of the $C_i$. If such a facet $F$ in not
available, $C$ is not covered by the $C_i$. The process must stop
because we only have finitely many affine spaces for splitting the
produced cones.

As for the group covering test, we first form the intersection $U$
of all the given full rank subgroups $G_1,\dots,G_m\subset\ZZ^r$.
Then we check whether an element of each the finitely many residue
classes in $\ZZ^r/U$ in $\ZZ^r$ belongs to one of the $G_j$.

Moreover, using the algorithm INTERSECTION in Subsection
\ref{VIRT} below, which computes intersections of modules with
affine subspaces, we can also give an algorithm for deciding
whether $M$ is asymptotically covered by $M_1,\dots M_t$ (again
using generating sets as input). One only needs to consider the
finite number of residue classes in $\ZZ^n$ modulo $\gp(S)$
represented by the given generators of $M$ -- their union contains
$M$.

The obtained algorithms, checking the asymptotic covering
condition both for semigroups and modules, will be called
ASYMPTOTIC.

We recall from \cite{BGcov} that the \emph{asymptotic
Carath\'eodory rank} $\CR^{\a}(S)$ of a positive affine semigroup
$S\subset\ZZ^n$ is defined as
$$
\min\bigl\{r:\ \pi(\{x\in S:\ \rho(x)\le r\}|S)=1\bigr\}.
$$
($\rho$ is the representation length, see Subsection
\ref{CANInt}), and the \emph{virtual Carath\'eodory rank}
$\CR^{\v}(S)$ is defined as
$$
\min\bigl\{r:\#(S\setminus\{x\in S:\ \rho(x)\le r\})<\infty\bigr\}.
$$

Lemma \ref{critas} has the following

\begin{corollary}\label{uh=ic}
\begin{itemize}
\item[(a)]
Suppose $S\subset\ZZ^n$ is an affine semigroup and $S_1,\dots,S_t$
are affine sub-se\-mi\-gro\-ups $S_1,\dots,S_t$ of $S$. If these
sub-semigroups are normal and cover $S$ asymptotically, then $S$
is normal and covered by $S_1,\dots,S_t$.
\item[(b)]
Assume $S\subset\ZZ^n$ is a positive affine semigroup. If
$\CR^{\a}(S)=\rank S$ then $S$ is normal, $\CR^{\v}(S)=
\CR(S)=\rank S$ and, moreover, $S$ satisfies (FHC). In particular,
(ICP) and (FHC) are equivalent and they imply the normality.
\item[(c)]
For $S$ as in upright there is an algorithm for computing $\CR(S)$
and, in particular, for checking (ICP) in finitely many steps.
\end{itemize}
\end{corollary}

\begin{proof}
Claim (a) is a direct consequence of Lemma \ref{critas}. Claim (b)
follows from the same lemma and the observation that if
$\CR^{\a}(S)=\rank S$, then the full rank free sub-semigroups of
$S$, generated by elements of $\Hilb(S)$, cover $S$
asymptotically. This is so because the contribution from
degenerate subsets of $\Hilb(S)$ is ``thin'' and cannot affect the
asymptotic covering property. (c) follows from (b) and
ASYMPTOTIC.
\end{proof}

\begin{remark}
A motivation for the introduction of asymptotic and virtual
Carath\'eodory ranks of positive semigroups is the following
improvement of Seb\H o's inequality \ref{SebCR}. Suppose $S$ is an
affine positive normal semigroup and $\rank S\geq3$; then
$$
\CR^{\a}(S)\leq2\rank S-3
$$
and if, in addition, $S$ is \emph{smooth}, then
$$
\CR^{\v}(S)\leq2\rank S-3.
$$
``Smooth'' here means $\ZZ x+S\approx\ZZ\dirsum\ZZ_+^{\rank S-1}$
for each extreme generator of $S$. (Equivalently, for a field $K$
the variety $\Spec(K[S])\setminus\{\frak m\}$ is smooth, where
$\frak m$ is the monomial maximal ideal of $K[S]$.) These
inequalities have been proved in \cite{BGcov}.
\end{remark}

\subsection{Virtual covers}\label{VIRT}

Now we develop an algorithm checking the virtual covering
condition. First we need an auxiliary algorithm that computes
intersections of semigroups and modules with affine spaces.

More precisely, assume $S\subset\ZZ^n$ is an affine semigroup and
$M\subset\ZZ^n$ is a finitely generated $S$-module, both given in
terms of generating sets, say $G_S$ and $G_M$. Let
$H_0\subset\RR^n$ be a rational subspace, given by a system of
rational linear forms, and $h\in\QQ^n$. By Gordan's lemma
$S_0=S\cap H_0$ is an affine semigroup and by \cite[7.2]{BGdiv}
$M_h=M\cap(h+H_0)$ is a finitely generated module over it. Our
goal is to find their generating sets.

By considering the intersections $(z+S)\cap(h+H_0)$, $z$ running
through $G_M$ one reduces the task to the special case when $M$ is
generated by a single element, i.~e. when $M$ is a parallel shift
of $S$ in $\ZZ^n$, say by $z$. Changing $M$ by $-z+M$ and $h$ by
$h-z$ we can additionally assume $M=S$. Furthermore, taking the
intersection $H_0\cap(\RR\tensor\gp(S))$, we may suppose that
$H_0\subset\RR\tensor\gp(S)$. In other words, it is enough to
consider the case $n=\rank S$.

Fix a surjective semigroup homomorphism $\phi:\ZZ^s_+\to S$,
$s=\#G_S$, mapping the standard generators of $\ZZ_+^n$ to the
elements of $G_S$. It gives rise to a surjective linear mapping
$\RR^s\to\RR^n$ which we denote again by $\phi$. Next we compute
$\Ker(\phi)$ and, using it, the preimage $L_0=\phi^{-1}(H_0)$ --
the latter is generated by $\Ker(\phi)$ and arbitrarily chosen
preimages of a basis of the rational space $H_0$. Then we find an
element $l\in\phi^{-1}(h)$. (Finding preimages requires only
solving linear systems of equations.)

When we have computed a generating set of the semigroup
$\ZZ^s_+\cap L_0$ and that of the module $\ZZ^s_+\cap(l+L_0)$ over
it, then, by applying $\phi$, we find the desired generating sets.
In other words, we have further reduced the problem to the special
case when $S=\ZZ^n_+$. The semigroup $\ZZ^n_+\cap H_0$ is normal
and positive. Its Hilbert basis is computed using NORMALIZ.

Next we check whether $\ZZ^n\cap(h+H_0)=\emptyset$. This is done
as follows. We compute a group basis $B_1$ of $H_0\cap\ZZ^n$ and
find a system of vectors $B_2$, disjoint form $B_1$, such that
$B_1\cup B_2$ is a basis of $\ZZ^n$. Then $B_2$ corresponds to a
basis of the real space $\RR^n/H_0$. We only need to check that the
residue class of $h$ is integral with respect to it. This is a
necessary and sufficient condition for
$\ZZ^n\cap(h+H_0)\not=\emptyset$.

If $\ZZ^n\cap(h+H_0)\neq\emptyset$, then we can pick a lattice
point $p$ in $\ZZ^n\cap(h+H_0)$. We declare it as the origin of
the affine subspace $p+H_0$ with the coordinate system represented
by $p+B_1$.

Next we compute the intersections $C=\RR^n_+\cap H_0$ and
$P=\RR^n_+\cap(h+H_0)$, representing them by systems of
inequalities in the coordinate systems of $H_0$ and $h+H_0$, which
are given by $B_1$ and $p+B_1$ respectively.

We can make the natural identification $H_0=\RR^m$, $m=\#B_1$.
Consider the convex hull $\Pi$ in $\RR^{m+1}$ of the subset
$$
(C,0)\cup(-p+P,1)\subset\RR^{m+1}.
$$
The crucial observation is that $\Pi$ is a finite rational pointed cone
(for a similar construction in the context of divisorial ideals
see \cite[Section 5]{BGdiv}). Then, using again NORMALIZ we compute
$\Hilb(\Pi)$. The last step consists of listing those elements of
$\Hilb(\Pi)$ which have 1 as the last coordinate. Returning to the
old copy of $\RR^n$ these elements represent the minimal
generating set of
$$
\RR^n_+\cap(h+H_0)\in\M(\RR^n_+\cap H_0).
$$
This algorithm will be called INTERSECTION.

Note that we do not exclude the case when $H_0\cap S=\{0\}$. Then
the algorithm above just lists the elements of the finite set
$M\cap(h+H_0)$.

Now assume $S_1,\dots,S_t\subset S$ and $M,M_1,\dots,M_t$ are as
in Subsection \ref{ASYMP}, given in terms of their generators. By
$\Sigma$, $C_\sigma$ and $z_i$ we refer to the same objects as in
Lemma \ref{critas}. We will describe an algorithm deciding the
virtual covering property for the given semigroups and modules. It
uses induction on $\rank S$.

In the case $\rank S=1$ one easily observes that asymptotic and
virtual covering conditions coincide by Lemma \ref{cond}. So we can apply 
ASYMPTOTIC.

Assume $\rank S>1$. Using ASYMPTOTIC we first check that
we have at least asymptotic covering.

Let us first consider the case of semigroups. For every
$\sigma\in\Sigma$ we can pick an element
$z_\sigma\in\bigcap_{i\in\sigma}(z_i+C(i))$. Then
$$
S_\sigma:=\gp(S)\cap(z_\sigma+C_\sigma)\subset S.
$$
An important observation is that the complement $(C_\sigma\cap
S)\setminus S_\sigma$ is contained in finitely many sets of the
type $(h+\RR F)\cap S$, where $h\in\gp(S)$ and $F\subset C_\sigma$
is a facet. ($\RR F$ refers to the linear space spanned by $F$.)
Moreover, we can list explicitly such affine subspaces
$h+\RR F$ that cover this complement. Namely, for any facet
$F\subset C_\sigma$ we consider a system of vectors
$$
\{h_0,h_1,\dots,h_{v_F(z_\sigma)}\}\subset\gp(S)
$$
satisfying the condition $v_F(h_j)=j$, $j\in[0,v_F(z_\sigma)]$,
where $v_F:\gp(S)\to\ZZ$ is the surjective group homomorphism
uniquely determined by the conditions $v(\RR F\cap\gp(S))=0$ and
$v_F((C_\sigma\cap\gp(S))\geq0$.

The semigroup $S$ is virtually covered by $S_1,\dots,S_t$ if and
only if $S\cap(h+\RR F)\in\M(S\cap\RR F)$ is virtually covered by
the modules
$$
S_i\cap(h+\RR F)\in\M(S_i\cap\RR F),\qquad i\in[1,t]
$$
for all the (finitely many) possibilities $\sigma\in\Sigma$,
$F\subset C_\sigma$ and $h\in\gp(S)$ as above.

All of these intersection semigroups and modules can be computed
with INTERSECTION. Therefore, having decreased the rank
by one, we can use induction.

In the case of modules we first reduce the general case to the
situation when $M\subset\gp(S)$ -- we just split the problem into
finitely many similar problems corresponding to the set of residue
classes of the given generators of $M$ modulo $\gp(S)$. We then
pick elements (say, among the given generators) $y_i\in M_i$,
$i\in[1,t]$ and also elements
$$
m_\sigma\in\bigcap_{i\in\sigma}(y_i+z_i+C(S_i)),\qquad\sigma\in\Sigma.
$$
We have
$$
M_\sigma:=\gp(S)\cap(m_\sigma+C_\sigma)\subset
M,\qquad\sigma\in\Sigma.
$$
Let $m$ be an element of the given generating set for $M$. Then
the complement $(M\cap(m+C_\sigma))\setminus M_\sigma$ is
contained in finitely many sets of the type $(h+\RR F)\cap M$,
where the $F$ are as above and the $h\in\gp(S)$ constitute a finite
system such that the $v_F(h)$ exhaust the integers between $v_F(m)$ and
$v_F(m_{\sigma})$. We see that all the steps we have
carried out for the semigroups can be performed in the situation
of modules as well -- we only need to go through the whole process
for every generator of $M$.

The produced algorithm, deciding the virtual covering property, is
called VIRTUAL.

\subsection{Covers}\label{COVERINGS}

Now we complete the algorithm deciding covering property for
semigroups and their modules, as mentioned at the beginning of
Section \ref{ALGOR}. The algorithm will be called COVERING.
Again, we use induction on rank of the big semigroup. Analyzing
VIRTUAL one observes that the inductive step in developing
COVERING can be copied word-by-word from VIRTUAL. So
the only thing we need to describe is COVERING for rank 1
semigroups.

Assume $S,S_1,\dots,S_t$ and $M,M_1,\dots,M_t$ are as above and,
in addition, $\rank S=1$. We restrict ourselves to the case when
$S$ is positive. The other case can be done similarly.

After computing $\gp(S)$, we can assume $\gp(S)=\ZZ$ without loss
of generality. Since $\ZZ$ is covered by a finite system of
subgroups exactly when one of the subgroups is the whole $\ZZ$ we
must check (according to Lemma \ref{critas})
that one of the groups $\gp(S_1),\dots,\gp(S_t)$ coincides with
$\ZZ$. Assume $\gp(S_1)=\ZZ$. By CONDUCTOR we find an element
$z\in\cc_{\bar S_1/S_1}$. Now we only need to make sure that
the finite set $[1,z]\cap S$ is in the union $S_1\cup\dots\cup S_t$.

For the modules we first reduce the general case to the situation
$M\subset\ZZ$ (as we did in the previous subsection) and, by a
suitable shift, further to the special case $0\in M\subset\ZZ_+$.
By Lemma \ref{critas} there is no loss of generality in assuming that
$M_1\not=\emptyset$. Then, again, we only have a finite problem
of checking that $[0,z+m]\cap M\subset M_1\cup\dots\cup M_t$,
where $z$ is as above and $m\in M_1$ is arbitrarily chosen element
(say, a given generator).

\begin{remark}
As mentioned, our goal in this section was to show that the
question whether a given affine semigroup is covered by a finite
system of affine sub-semigroups can be checked algorithmically.
However, we did not try to make the algorithm as optimal as
possible. For instance, our arguments
use heavily conductor ideals and we work with random elements in
these ideals. On the other hand in some special cases one can
compute $\cc_{\bar S/S}$ exactly. Especially this is possible in
the situation when $S$ is a positive affine semigroup, generated
by $\rank S+1$ elements; see \cite{RR}.

The real motivation for implementing a part of the algorithms
above would be a semigroup that violates (UHC), but resists all
random tests for detecting the violation of (ICP) (or,
equivalently, (FHC); see Corollary \ref{uh=ic}(b)). Unfortunately,
so far we have only found 2 essentially different semigroups
violating (UHC), and they violate (ICP) too.
\end{remark}

\section{Algebraic properties of affine semigroup algebras} \label{Trung}

In this section we always consider affine semigroups $S$ of
$\ZZ_+^r$ (often $r$ will be the rank of $S$, but we do not
necessarily assume this). Then the affine semigroup algebra $K[S]$
over a field $K$ can be viewed as a subalgebra of the polynomial
ring $K[T_1,\dots,T_r]$.

\subsection{Defining equations}

Let $\Hilb(S) = \{x_1,\ldots,x_n \}$. Consider the semigroup
homomorphism $\pi: \Bbb Z_+^n \to S$ given by $(u_1,\ldots,u_n)
\mapsto u_1x_1 + \ldots + u_nx_n$. Let $K[X] = K[X_1,\dots,X_n]$
be a polynomial ring over a field $K$ in $n$ indeterminates. The
map $\pi$ lifts to a homomorphism of semigroup algebras $\phi:
K[X] \to K[S]$. The kernel of $\phi$ is a prime ideal $I_S$ in
$K[X]$ and we have a representation of the semigroup algebra
$$
K[S] \cong K[X]/I_S.
$$
The ideal $I_S$ is often called the \emph{toric ideal} of $S$. The
following result is well-known (for example, see Gilmer
\cite{Gi}).

\begin{proposition}
The toric ideal $I_S$ is generated by the set of binomials
$$
\{X^u-X^v|\ u,v \in \Bbb Z_+^n\ \text{\rm with}\ \pi(u) =
\pi(v)\}.
$$
\end{proposition}

Let $\mu(I)$ denote the minimal number of generators of an ideal
$I$. Because of the above property of $I_S$ one might think that
$\mu(I_S)$ could not be big or, more precisely, that $\mu(I_S)$
were bounded by a number which depends only on the number $n$. But
that is not the case.

Let $S$ be a \emph{numerical semigroup}, that is $S \subseteq \Bbb
Z_+$. If $n = 2$, then $\mu(I_S) = 1$ because $I_S$ is a principal
ideal in $K[X_1,X_2]$. If $n = 3$, Herzog \cite{He} proved that
$\mu(I_S) \le 3$. If $n \ge 4$, Bresinsky \cite{Bre1} showed that
$\mu(I_S)$ can be arbitrarily large.

However, one may expect that $\mu(I_S)$ depends only on $n$ for
special classes of affine semigroups. Let $S$ be generated by
$n$ non-negative integers $x_1,\ldots,x_n$. Without restriction we
may assume that $x_1,\ldots,x_n$ have no common divisor other than
1. Then there exists an integer $c$ such that $a \in S$ for all
integers $a \ge c$ (i.~e. $c$ is in the conductor ideal). Let $c$
be the least integer with this property. We call $S$ a \emph{symmetric
numerical semigroup} if $a\in S$ whenever $c-a-1 \not\in S$,
$a \in \Bbb N$.

\begin{example}
Let $S = \langle 6,7,8\rangle$. Then
$$
S = \{0,6,7,8,12,13,14,15,16,18,19,20,21,22,\dots\}.
$$
Hence $c = 18$. It is easy to check that $S$ is a symmetric numerical semigroup.
\end{example}

The interest on symmetric numerical semigroups originated from the
classification of plane algebroid branches \cite{Ap}. Later, Herzog and
Kunz \cite{HK} realized that symmetric numerical semigroups correspond
to Gorenstein affine monomial curves.

\begin{problem}
Let $S$ be a symmetric numerical semigroup. Does there exist an
upper bound for $\mu(I_S)$ which depends only on the minimal
number of generators of $S$?
\end{problem}

If $n = 3$, Herzog \cite{He} proved that $\mu(I_S) = 2$. If $n =
4$, Bresinsky \cite{Bre2} proved that $\mu(I_S) \le 5$. If $n =
5$, Bresinsky \cite[Theorem 1]{Bre3} proved that $\mu(I_S) \le
13$, provided $x_1+x_2 = x_3+x_4$. It was also Bresinsky
\cite[p.~218]{Bre2}, who raised the above problem which has
remained open until today.

Instead of estimating the number of generators of $I_S$ one can
also try to bound the degree of the generators. We will discuss
this problem in Subsections \ref{Koszul} and  \ref{CastMum}.

We call the least integer $s$ for which there exist binomials
$f_1,\ldots,f_s$ such that $I_S$ is the radical of the ideal
$(f_1,\ldots,f_s)$ the \emph{binomial arithmetical rank} of $I_S$
and we will denote it by $\bara(I_S)$. Geometrically, this means
that the affine variety defined by $I_S$ is the intersection of
the hypersurfaces $f_1 = 0,\dots,f_s = 0$. In general, we have
$\hht I_S \le \bara(I_S) \le \mu(I_S)$.

\begin{problem}
Does there exist an upper bound for $\bara(I_S)$ in terms of $n$?
\end{problem}

We mention only a few works on this problem. If $S$ is a
homogeneous (i.e.\ graded and generated by elements of degree $1$)
affine semigroup in $\Bbb Z_+^2$, Moh \cite{M} proved that
$\bara(I_S) = n-2$ for $K$ of positive characteristic. This
implies that $I_S$ is a set-theoretic complete intersection. If
$K$ has characteristic $0$ and $S$ is as above, then Thoma
\cite{Th} has shown that $\bara(I_S)=n-2$ if $I_S$ is a complete
intersection, otherwise $\bara(I_S)=n-1$. These results have been recently
generalized by Barile, Morales and Thoma [BMT] to affine semigroup algebras of the form
$$K[S] = K[t_1^{d_1},...,t_r^{d_r},
t_1^{a_{11}} \cdots t_r^{a_{1r}},...,t_1^{a_{s1}} \cdots t_r^{a_{sr}}],$$
where $d_1,...,d_r$ and $a_{11},...,a_{sr}$ are positive integers.

\subsection{Initial ideals and the Koszul property}\label{Koszul}

Let $K[X] = K[X_1,\dots,X_n]$ be a polynomial ring over a field
$k$. As usual, we will identify a monomial $X^u = X_1^{u_1}\cdots X_n^{u_n}$
with the lattice point $u = (u_1,\ldots,u_n)$. A total order $<$
on $\ZZ_+^n$ is a \emph{term order} if it has the following
properties:
\begin{itemize}
\item[(i)] the zero vector 0 is the unique minimal
element;
\item[(ii)] $v < w$ implies $u+v < u+w$ for all $u,v,w \in
\Bbb N^n$.
\end{itemize}
Given a term order $<$, every non-zero polynomial $f \in K[X]$ has
a largest monomial which is called the \emph{initial monomial} of
$f$. If $I$ is an ideal in $K[X]$, we denote by $\ini(I)$ the
ideal generated by the initial monomials of the elements of $I$.
This ideal is called the \emph{initial ideal} of $I$. The passage
from $I$ to $\ini(I)$ is a flat deformation (see e.g.
\cite[15.8]{Ei}). Hence one can study $I$ be means of $\ini(I)$.

With every monomial ideal $J$ we can associate the following
combinatorial object
$$
\Delta(J) := \{F \subseteq \{1,\dots,n\}:\ \text{there is no
monomial in $J$ whose support is $F$}\}
$$
where the support of a monomial $X^a$ is the set $\{i:a_i\neq
0\}$. Clearly $\Delta(J)$ is a simplicial complex on the vertex
set $\{1,\dots,n\}$, and it easily seen that $J$ and its radical
$\sqrt{J}$ define the same simplicial complex: $\sqrt{J}$ is
generated by all \emph{square-free monomials} $X_{i_1}\cdots
X_{i_s}$, $i_1<\dots< i_s$, for which $\{i_1,\dots,i_s\}$ is not a
face of $\Delta(J)$.

We call $\Delta(\ini(I))$ the \emph{initial complex} of $I$ (with
respect to the term order $<$).

For a toric ideal $I_S$ one may ask whether there is a
combinatorial description of the initial ideal $\ini(I_S)$ or, at
least, their radicals.

In the remaining part of this subsection we assume that $S$ is a
homogeneous affine semigroup $S_M\subset\ZZ^{r+1}$ with Hilbert basis
$M=\{x_1,\dots,x_n\}\subset \ZZ^r$. By $C$ we denote the cone
$C(S)$.

An $M$-triangulation
of $C$ is called \emph{regular} if there is a \emph{weight vector}
$\omega = (\omega_1,\dots,\omega_n) \in \RR^n_+$ such that the
simplicial cones of the triangulation are spanned exactly by those
subsets $F\subset M$ for which there exists a vector $c \in \Bbb
R^r$ with
\begin{align*}
\langle c,x_i\rangle = \omega_i & \quad \text{ if }  x_i \in F,\\
\langle c,x_j\rangle < \omega_j & \quad \text{ if}  x_j \not\in F.
\end{align*}
Geometrically, the simplicial cones of a regular triangulation of
$C(S)$ are the projections of the lower faces of the convex hull
$P$ of the vectors
$\bigl\{(x_1,\omega_1),\dots,(x_n,\omega_n)\bigr\}$ in $\Bbb
R^{r+1}$ onto the first $r$ coordinates. Note that a face of $P$
is lower if it has a normal vector with negative last coordinate.

It is clear that every $M$-triangulation of $C(S)$ can be
identified with the simplicial complex of those subsets
$\{i_1,\dots,i_r\} \subseteq \{1,\dots,n\}$ for which the vectors
$x_{i_1},\ldots,x_{i_r}$ span a face of a simplicial cone of the
triangulation. Using this identification, Sturmfels \cite{Stu},
Theorem 8.3 and Corollary 8.4, discovered the following
connections between the initial complexes of $I_S$ and the
triangulations of $C(S)$.

\begin{theorem}
The initial complexes $\Delta(\ini(I_S))$ are exactly the
simplicial complexes of the regular $M$-triangulations of
$C(S)$.
\end{theorem}

\begin{corollary}
The ideal $\ini(I_S)$ is generated by square-free monomials if and
only if the corresponding regular $M$-triangulation of $C(S)$ is
unimodular.
\end{corollary}

Therefore, if $I_S$ has a square-free initial ideal, then $S$ must be
normal and, being generated in degree $1$, polytopal (see Proposition
\ref{UHCsat}). On the other hand, as the counterexample in Subsection
\ref{CANInt} shows, there exist normal lattice polytopes without any
unimodular triangulation (even without unimodular covering). Therefore
$I_S$ need not have a square-free initial ideal for normal polytopal
semigroups $S$. (There also exist polytopes that have a unimodular
triangulation, but no such regular triangulation; see Ohsugi and Hibi
\cite{OH1}.)

However, as observed in Subsection \ref{CANInt}, any triangulation
of a lattice polytope of dimension 2 into empty lattice simplices
is unimodular by Pick's theorem, and $I_S$ has plenty of
square-free initial ideals in this special situation.

The results of Sturmfels give us a method to prove that a
semigroup algebra is Koszul. Recall that a homogeneous algebra $A$
over a field $K$ is called \emph{Koszul} if $K$ as an $A$-module
has a resolution:
$$
\cdots \longrightarrow E_2 \overset {\phi_2}
\longrightarrow E_1 \overset {\phi_1} \longrightarrow A \longrightarrow K
\longrightarrow 0,
$$
where $E_1, E_2, \dots$ are free $R$-modules and the entries of
the matrices $\phi_1, \phi_2,\dots$ are forms of degree 1 in $A$.
For more information see the survey of Fr\"oberg \cite{Fr}.

Let $A = R/I$ be a presentation of $A$, where $R$ is a polynomial
ring over $K$ and $I$ is a homogeneous ideal in $R$. If $A$ is a
Koszul algebra, then $I$ must be generated by quadratic forms. The
converse is not true. However, $A$ is Koszul if there exists a
term order $<$ such that the initial ideal $\ini(I)$ is generated
by quadratic monomials. Therefore, if a lattice
polytope $P$ has a unimodular regular triangulation whose minimal
non-faces are edges, then the semigroup algebra $K[P]$ is Koszul.

We have proved in \cite{BGT} that the following classes of lattice
polytopes have this property:
\begin{itemize}
\item[(1)] lattice polytopes in $\Bbb
R^2$ whose boundaries have more than 3 lattice points,
\item[(2)] lattice polytopes in $\Bbb R^r$ whose facets are parallel to the
hyperplanes given by the equations $T_i = 0$ and $T_i-T_j = 0$.
\end{itemize}
In particular, it can be shown that if $P$ is a lattice polytope
in $\Bbb R^2$ with more than 3 lattice points, then $K[P]$ is
Koszul if and only if the boundary of $P$ has more than 3 lattice
points.

It would be of interest to find more lattice polytopes which have
unimodular regular triangulations whose minimal non-faces are
edges. For any lattice polytope $P \subset \Bbb R^r$, it is known
that the semigroup algebra $K[cP]$ is Koszul for $c \ge r$
\cite[Theorem 1.3.3]{BGT}. This has led us to the following
problem.

\begin{problem}
Does $cP$, $c \gg 0$, have a unimodular regular triangulation
$\Delta$ such that the minimal non-faces of $\Delta$ are edges?
\end{problem}

We have already stated this problem in Subsection \ref{HighMult},
however without the attribute ``regular'' and the condition that
the minimal non-faces of $\Delta$ should be edges. In this
connection we have pointed out that unimodular triangulations for
$cP$ have been constructed for infinitely many $c$ in \cite{KKMS};
these triangulations are in fact regular.

It has been asked whether a Koszul semigroup algebra always has an
initial ideal generated by quadratic monomials. But this question
has a negative answer by Roos and Sturmfels \cite{RS}. There also exist
normal non-Koszul semigroup algebras defined by quadratic binomials; see
Ohsugi and Hibi \cite{OH2}.

\subsection{The Cohen-Macaulay and Buchsbaum properties}

Let $(A,\mm)$ be a local ring. A system of elements
$x_1,\dots,x_s$ of $A$ is called a \emph{regular sequence} if
$$
(x_1,\ldots,x_{i-1}): x_i = (x_1,\dots,x_{i-1}),\qquad\ i =
1,\dots,s.
$$
It is called a \emph{weak-regular sequence} if
$$
\mm\bigl[(x_1,\ldots,x_{i-1}): x_i\bigr] \subseteq
(x_1,\dots,x_{i-1}),\qquad i = 1,\dots,s.
$$
Let $d = \dim A$. A system of $d$ elements $x_1,\dots,x_d$ of $A$
is called a \emph{system of parameters} of $A$ if the ideal
$(x_1,\dots,x_d)$ is an $\mm$-primary ideal. The local ring $A$ is
called a \emph{Cohen-Macaulay ring} if there exists an (or every)
system of parameters of $A$ is a regular sequence. It is called a
\emph{Buchsbaum ring} if every system of parameters of $A$ is a
weak-regular sequence. If $A$ is a finitely generated homogeneous
algebra over a field and $\mm$ is its maximal homogeneous ideal,
then we call $A$ a Cohen-Macaulay resp.\ Buchsbaum ring if the
local ring of $A$ at $\mm$ is  Cohen-Macaulay resp. Buchsbaum.
Cohen-Macaulay resp. Buchsbaum rings can be characterized in
different ways and they have been main research topics in
Commutative Algebra. See \cite{BH} and \cite{SV} for more
information on these classes of rings.

By a fundamental theorem of Hochster \cite{Ho} normal affine
semigroup rings are Cohen-Macaulay. For general affine semigroup
rings the Cohen-Macaulay property has been characterized in
\cite[Theorem 3.1]{TH}, which is based on earlier work of Goto and
Watanabe \cite{GW}. For two subsets $E$ and $F$ of $\Bbb Z^r$ we
set
$$
E \pm F = \{v \pm w|\  v \in E, w \in F\}.
$$
Let $F_1,\dots,F_m$ be the facets of the cone $C(S)$. Put $S_i = S
- (S\cap F_i)$ and
$$
S' = \bigcap_{i = 1}^mS_i.
$$
For every subset $J$ of the set $[1,m] = \{1,\dots,m\}$ we set
$$
G_J = \bigcap_{i \not\in J}S_i \setminus \bigcup_{j\in J}S_j,
$$
and we denote by $\pi_J$ the simplicial complex of non-empty
subsets $I$ of $J$ with $\bigcap_{i \in I}(S \cap F_i) \neq
\{0\}$.

\begin{theorem}
$K[S]$ is a Cohen-Macaulay ring if and only if the
following conditions are satisfied:
\begin{itemize}
\item[(a)] $S' = S$;
\item[(b)] $G_J$ is either empty or acyclic over $K$
for every proper subset $J$ of $[1,m]$.
\end{itemize}
\end{theorem}

Though Buchsbaum rings enjoy many similar properties like those of
Cohen-Macaulay rings, one has been unable to find a similar
characterization for the Buchsbaum property of $K[S]$.

\begin{problem}
Find criteria for an affine semigroup algebra $K[S]$ to be a
Buchsbaum ring in terms of the affine
semigroup $S$.
\end{problem}

Recall that an affine semigroup $S$ is called \emph{simplicial} if
$C(S)$ is spanned by $r$ vectors of $S$, where $r = \rank S$. Geometrically,
this means that $C(S)$ has $r$ extreme rays or, equivalently, $r$
facets. This class contains all affine semigroups in $\Bbb Z^2$.
Goto, Suzuki and Watanabe \cite{GSW} resp.\ Trung \cite{Tr1} gave the following simple criteria for a simplicial affine semigroup algebra to be Cohen-Macaulay resp.\ Buchsbaum.

\begin{theorem}
Let $S$ be a simplicial affine semigroup with $d = \rank \gp(S)$.
Let $v_1,\ldots,v_d$ be the vectors of $S$ which span $C(S)$. Then
\begin{itemize}
\item[(a)] $K[S]$ is Cohen-Macaulay if and only if
$$
\{v \in \gp(S): v+v_i, v+v_j \in S\ \text{\rm for some indices
}i\neq j\} =S;
$$
\item[(b)] $K[S]$ is Buchsbaum if and only if
$$
\{v \in \gp(S)|\ v+2v_i, v+2v_j \in S\
\text{\rm for some indices $i \neq j$}\} + \Hilb(S) \subseteq S.
$$
\end{itemize}
\end{theorem}

The above criteria are even effective. For example consider (a).
Then we form the intersection
$$
(-s_i+S)\cap (-s_j+S)
$$
of $S$-modules and test whether this module is contained in $S$.
Section \ref{ALGOR} contains algorithms for these tasks. From the
ring-theoretic point of view, the main special property of
simplicial affine semigroups is the existence of a homogeneous
system of parameters consisting of monomials. Therefore certain
homological properties that depend on system of parameters can be
formulated in terms of the semigroup.

 What we know on a given affine semigroup is usually its
Hilbert basis. Therefore, we raise the following stronger problem.

\begin{problem} Find criteria for $K[S]$ to be a Cohen-Macaulay or
Buchsbaum ring in terms of $\Hilb(S)$.
\end{problem}

This problem is not even solved for the class of homogeneous
affine semigroups in $\Bbb Z_+^2$ which are generated by subsets
of
$$
M_e = \{v = (v_1,\dots,v_r) \in \Bbb Z_+^r|\ v_1 + \cdots + v_r
= e\},
$$
where $e$ is a given positive number. The algebra of the semigroup
generated by the full set $M_e$ is just the homogeneous coordinate
ring of the $e$-th Veronese embedding of the
$(r-1)$-dimensional projective space. The algebras generated by
subsets of $M_e$ are the homogeneous coordinate rings of projections of
this Veronese variety.

Gr\"obner \cite{Gr} was the first who studied the Cohen-Macaulay
property of such semigroup algebras. Let $H$ be an arbitrary
subset of $M_e$ and $S = \langle H \rangle$. If $H$ is obtained from
$M_e$ by deleting one, two, or three vectors, we know exactly when
$K[S]$ is a Cohen-Macaulay or Buchsbaum ring \cite{Sc, Tr1, Hoa}.
If $r = 2$, we may identify $H$ with the sequence
$\alpha_1,\dots,\alpha_n$ of the first coordinates of the vectors
of $H$. There have been some attempts to determine when $K[S]$ is
a Buchsbaum or Cohen-Macaulay ring in terms of
$\alpha_1,\dots,\alpha_n$. But satisfactory answers were obtained
only in a few special cases \cite{Bre4, BSV, Tr2}.

\subsection{Castelnuovo-Mumford regularity}\label{CastMum}

Let $A = \bigoplus_{t\ge 0}A_t$ be a finitely generated homogeneous
algebra over the field $K$. Let $A = R/I$ be a representation of
$A$, where $R$ is a polynomial ring over $K$ and $I$ a
homogeneous ideal of $R$. Then we have a finite minimal free
resolution of $A$ as a graded $R$-module:
$$
0 \longrightarrow E_s \longrightarrow \cdots \longrightarrow E_1
\longrightarrow R \longrightarrow A \longrightarrow 0,
$$
where $E_1,\dots,E_s$ are graded $R$-modules. Let $b_i$ be the
maximum degree of the generators of $E_i$, $i = 1,\dots,s$. Then
the {\it Castelnuovo-Mumford regularity} of $A$ is defined as the
number
$$
\reg(A) := \max\{b_i-i|\ i = 1,\dots,s\}.
$$
It is independent of the representation of $A$. In fact, it can be
defined solely in terms of $A$ as follows.

Let $\mm$ denote the maximal homogeneous ideal of $A$. For any
$A$-module $M$ we set
$$
\Gamma_\mm(M) := \{x \in M|\ x\mm^t = 0\ \text{\rm for some number $t \ge 0$}\}.
$$
Then $\Gamma_\mm(*)$ is a left exact additive functor from the
category of $A$-modules into itself. Let $H_\mm^i(*)$ denote the
$i$-th right derived functor of $\Gamma_\mm(*)$. Then $H_\mm^i(M)$
is called the $i$-th \emph{local cohomology module} of $M$ (with
respect to $\mm$). If $M$ is a graded $A$-module, then
$H_\mm^i(M)$ is also a graded $A$-module. Write $H_\mm^i(M) =
\bigoplus_{t \in \Bbb Z}H_\mm^i(M)_t$. It is known that  $\reg(A)$ is
the least integer $m$ such that $H_\mm^i(A)_t = 0$ for all $t >
m-i$ and $i \ge 0$.

The Castelnuovo-Mumford regularity $\reg(A)$ is an extremely
important invariant because it is a measure for the complexity of
$A$. For instance, $\reg(A)+1$ is an upper bound for the maximal
degree of the defining equations of the ideal $I$. See \cite{EG} and
\cite{Ei} for more information on the Castelnuovo-Mumford regularity.

It is a standard fact that the \emph{Hilbert function} $\dim_K A_t$
is a polynomial $P_A(t)$ of degree $d-1$ for $t \gg 0$, where $d =
\dim A$. If we write
$$
P_A(t) = \frac{et^{d-1}}{(d-1)!} +\ \text{\rm terms of degree $< d-1$},
$$
then $e$ is called the \emph{multiplicity} of $A$. Let $n$ denote
the minimal number of generators of $A$. In general, the
regularity is bounded by a double exponential function of $e$, $d$
and $n$. However, if $A$ is a domain, there should be an upper
bound for $\reg(A)$ with lower complexity. In this case, Eisenbud
and Goto \cite{EG} have conjectured that
$$
\reg(A) \le e + n - d.
$$
Gruson, Lazarsfeld and Peskine proved this conjecture in the case
$\dim A = 2$ \cite{GLP} (the case $\dim A = 1$ is trivial). For
$\dim A \ge 3$, it has been settled only under some additional
conditions on $A$.

For affine semigroup algebras, the above conjecture is still open.
Except those cases which can be derived from the known results for
homogeneous domains, the conjecture has been settled only for
affine semigroup algebras of codimension 2 ($n-d = 2$) by Peeva
and Sturmfels \cite{PS}.

Let $S$ be a homogeneous affine semigroup. Then $\Hilb(S)$ must
lie on a hyperplane of $\Bbb R^r$. It is known that the
multiplicity $e$ of $K[S]$ is equal to the normalized volume of
the convex polytope spanned by $\Hilb(S)$ in this hyperplane (for
example, see \cite[6.3.12]{BH}). Moreover, one can also describe
the regularity \cite{Stu} and the local cohomology of $K[S]$
combinatorially in terms of $S$ (see e.g. \cite{St2, TH, SS}
or \cite[Ch.\ 6]{BH}).

For a homogeneous affine semigroup $S$ in $\Bbb
Z_+^2$, this description is very simple. Without restriction we
may assume that $\Hilb(S)$ consists of vectors of the forms
$(a,1)$, $0 \le a \le e$, where
the vectors $v_1 = (0,1)$ and $v_2 = (e,1)$ belong to
$S$. Let $S'$ denote the set of vectors $v \in \Bbb Z_+^2$ for
which there are positive integers $m_1,m_2$ such that $v + m_1v_1
\in S$, $v + m_2v_2 \in S$. Then $\reg(K[S]) = \max\{a+b|\ (a,b)
\in S'\setminus S\} + 1$. By the result of Gruson, Lazarsfeld and Peskine,
$\reg(K[S]) \le e - n+2$, where $n$ is the number of vectors of
$\Hilb(S)$. It would be nice if we could find a combinatorial
proof for this bound.

We say that a binomial $X^{u}-X^{v} \in I_S$ is \emph{primitive}
if there is no other binomial $X^{u'} - X^{v'} \in I_S$ such that
$X^{u'}$ divides $X^{u}$ and $X^{v'}$ divides $X^{v}$. The set of
all primitive binomials of $I_S$ generates $I_S$. It is called the
\emph{Graver basis} of $I_S$ and denoted by $\Gr_S$. A binomial
$\zeta$ in $I_S$ is called a \emph{circuit} of $S$ if its support
$\supp(\zeta)$ (the set of variables appearing in $\zeta$) is
minimal with respect to inclusion. The \emph{index} of a circuit
$\zeta$ is the index of the additive group generated by
$\supp(\zeta)$ in the intersection of $\gp(S)$ with the linear
space spanned by $\supp(\zeta)$ in $\Bbb R^n$.

\begin{problem}
Prove that the degree of every binomial in $\Gr_S$ is bounded
above by the maximum of the products of the degree and the index
of the circuits of $S$.
\end{problem}

This problem was raised by Sturmfels in \cite[Section 4]{Stu2}. If it
has a positive answer, then one can show that $K[S]$ is defined by
binomials of degree $\le e-1$.

Suppose now that $P$ is a normal lattice polytope of dimension
$d$. Then the Castelnuovo-Mumford regularity of $R=K[P]$ has a
very simple geometric description. In fact,
$$
\reg(R)=d+1-\ell
$$
where $\ell$ is the minimal degree of a lattice point in the
interior of $C(P)$. In particular, one always has $\reg(A)\le d$,
and it follows that the ideal $I=I_{S_P}$ is generated by
binomials of degree $\le d+1$. It easily seen that the bound $d+1$ is attained if $P$ is a simplex (i.e.\ spanned by $d+1$ lattice points) with at least one lattice point in its interior, but no
lattice points in its boundary different from its vertices.
However, no counterexample seems to be known to the following
question:

\begin{problem}
Let $P$ a normal lattice polytope of dimension $d$ whose boundary
contains at least $d+2$ lattice points. Is $K[P]$ defined by
binomials of degree $\le d$?
\end{problem}

Clearly the answer is ``yes'' if $P$ has no interior lattice
point, and as we have seen in the previous section, it is also
``yes'' for $d=2$ since for $d=2$ one can even find a Gr\"obner
basis of $I$ of binomials of degree $2$ if $P$ contains at least
$4$ lattice points in its boundary. As far as the combinatorics of
triangulations is concerned, the result can be extended to higher
dimension. In fact, one has the following theorem
(\cite[3.3.1]{BGT})

\begin{theorem}
Let $P$ be a lattice polytope of dimension $d$ with at least $d+2$
lattice points in its boundary and at least one interior lattice
point. Then $P$ has a regular triangulation $\Delta$ into empty
lattice simplices such that the minimal non-faces of $\Delta$ have
dimension $\le d-1$.
\end{theorem}

Since one cannot expect the triangulation to be
unimodular for $d \ge 3$, the theorem only bounds the degree of the
generators of $\sqrt{\ini(I)}$. Nevertheless, one should
strengthen the last problem as follows:

\begin{problem}
Let $P$ a normal lattice polytope of dimension $d$ whose boundary
contains at least $d+2$ lattice points. Does $I_{S_P}$ have a
Gr\"obner basis consisting of binomials of degree $\le d$?
\end{problem}

We would like to mention that Sturmfels already raised in [Stu2] the
conjecture that for any normal lattice polytope of dimension $d$,  there exists
a Gr\"obner basis for $I_{S_P}$ consisting of binomials  of degree $\le d+1$.
There one can find some interesting problems on the maximal degree of the
defining equations and the regularity of toric ideals.

\end{document}